\documentclass{amsart}
\usepackage{radon}

\begin{document}

\title{Extending the Support Theorem to Infinite Dimensions}

\author[J. J. Becnel]{Jeremy J.~Becnel}
\address{Department of Mathematics and Statistics, Stephen F. Austin State
University, Nacogdoches, Texas 75962-3040}
\email{becneljj@sfasu.edu}
  

\date{\today}

\begin{abstract} The Radon transform is one of the most useful and applicable tools in functional analysis. First constructed by John Radon in 1917 \cite{JR17} it has now been adapted to several settings. One of the principle theorems involving the Radon transform is the Support Theorem. In this paper, we discuss how the Radon transform can be constructed in the white noise setting. We also develop a Support Theorem in this setting.
\end{abstract}
\maketitle


 
\section{Introduction} 
 
 The Radon transform \cite{JR17} associates to a function $f$ on the finite-dimensional space $\Rn$  a function
$R_f$ on the set of all hyperplanes in $\Rn$ whose value on any hyperplane $P$ is the integral of $f$ over $P$:
\begin{equation}\label{defRf}
Rf(P)=\int_Pf(x)\,dx,\end{equation}
 the integration here being with respect to Lebesgue measure on $P$. This transform does not generalize directly to infinite dimensions
because there is no useful notion of Lebesgue measure in infinite dimensions. However, there is a well-developed theory
of Gaussian measures in infinite dimensions and so it is natural to extend the Radon transform to infinite dimensions using
Gaussian measure:
\begin{equation}
Gf(P)=\int f\,d\mu_P,\end{equation}
where $\mu_P$ is the Gaussian measure on any infinite dimensional hyperplane $P$ in a Hilbert space $H_0$. A version of this transform was developed
in \cite{vM00} but we shall present a another account below.  A central feature of the classical
Radon transform $R$ is the Support Theorem (see, for instance, Helgason \cite{SH99}):
\begin{thm}[Support Theorem] \label{T:ST}
If $f$ is a rapidly decreasing continuous function
for which $R_f(P)$ is $0$ on every hyperplane $P$ disjoint from some compact convex set $K$  then $f(x)=0$ for $x \notin K$.
\end{thm}
\noindent or more appropriately for the Gaussian measure we have
\begin{thm}[Support Theorem---Gaussian] \label{T:STGaus}
If $f$ is a exponentially bounded continuous function
for which $G_f(P)$ is $0$ on every hyperplane $P$ disjoint from some compact convex set $K$  then $f(x)=0$ for $x \notin K$.
\end{thm}
 In this paper
we prove the infinite-dimensional version of this Support Theorem.

\section{White Noise Setup} \label{ss:taf}

We begin by describing the setting under which White Noise Analysis takes place. Because of the absence of the Lebesgue measure in infinite dimensions we begin by constructing a Gaussian measure. From here we can develop appropriate sets of tests functions and distributions.

 We work with   a
real separable Hilbert space $H_0$, and  a positive Hilbert-Schmidt
operator $A$ on $H_0$ such that there is orthonormal basis
$\{e_n\}_{n=1}^{\infty}$ of eigenvectors of $A$ and eigenvalues $\{\gl_n\}_{n=1}^{\infty}$ satisfying  
\begin{enumerate}
\item $Ae_n=\lambda_ne_n$
\item $1 < \gl_1 < \gl_2 < \dots$
\item $\sum_{n=1}^{\infty} \gl_n^{-2} < \infty$
\end{enumerate}
 
The example to keep in mind is
\begin{eqnarray*}
H_0&=& L^2({\mathbb{R}})\\
e_n&=& \phi_n\\
 A&=&-\frac{d^2}{dx^2}+\frac{x^2}{4}+\frac{1}{2}\quad\mbox{
 with
eigenvalues $\lambda_n = (n+1)$. }\end{eqnarray*}

We have the coordinate map
$$J: H_0\mapsto {\mathbb{R}}^W: f\mapsto \bigl(\langle
f,e_n\rangle\bigr)_{n\in W}$$
where we use the notation $W=\{1,2,...\}$. Let
\begin{equation}\label{eq:defFzero}
F_0=J(H_0)=\{(x_n)_{n\in W}: \sum_{n\in
W}x_n^2<\infty\}\end{equation} Now, for each $p\in W$, let
\begin{equation}\label{eq:defFp}F_p=\{(x_n)_{n\in W}:
\sum_{n\in W}\lambda_n^{2p}x_n^2<\infty\}
\end{equation}
On $F_p$ we have the inner-product $\langle\cdot,\cdot\rangle_p$
given by
$$\langle a, b\rangle_p =\sum_{n\in W}\lambda_n^{2p}a_nb_n$$
This makes $F_p$ a real Hilbert space, unitarily isomorphic to
$L^2(W, \mu_p)$ where $\mu_p$ is the measure on $W$ specified by
$\mu_p(\{ n\})=\lambda_n^{2p}$. Moreover, we have
\begin{equation}\label{eq:chainFp}
F\stackrel{\rm def}{=}\bigcap_{p\in W}F_p\subset\cdots F_2\subset
F_1\subset F_0=L^2(W,\mu_0)
\end{equation}
and each inclusion $F_{p+1}\to F_p$ is Hilbert-Schmidt.

Now we pull all this back to $H_0$. First set
\begin{equation}\label{eq:defEpnew}
H_p=J^{-1}(F_p)=\{x\in H_0: \sum_{n \in W}\lambda_n^{2p}|\langle
x,e_n\rangle|^2<\infty\} 
\end{equation}
 It is readily checked that
$H_p=\{x\in H_0 \, ; \, \norm|x|_p < \infty \}$ where $\norm|x|_p = |Ax|_0$ and also $H_p = A^{-p}(H_0)$. On $H_p$ we have the pull back
inner-product $\langle\cdot,\cdot\rangle_p$, which works out to
\begin{equation}\label{eq:Epinnerprod}
\langle f,g\rangle_p=\langle A^{p}f, A^{p}g\rangle
\end{equation} Then we have the chain
\begin{equation}\label{eq:chainEp}
\Hnucl \stackrel{\rm def}{=}\bigcap_{p\in W}H_p\subset\cdots H_2\subset
H_1\subset H_0,
\end{equation}
with each inclusion $H_{p+1}\to H_p$ being Hilbert-Schmidt.

Equip $\Hnucl$    with the topology generated by the norms
$|\cdot|_p$ (i.e. the smallest topology making all inclusions
$\Hnucl \to H_p$ continuous). Then $\Hnucl$  is, more or less by definition, a
nuclear space.
 The vectors $e_n$ all lie in $\Hnucl$  and the set of all
 rational-linear combinations of these vectors produces a
 countable dense subspace of $\Hnucl$.

 Consider a linear functional on $\Hnucl$  which is continuous. Then it
 must be continuous with respect to some norm $|\cdot|_p$.
 Thus
 the topological dual $\Hstar$  is the union of the duals $H'_p$. In
 fact, we have:
\begin{equation}\label{eq:chainEpdual}
\Hstar=\bigcup_{p\in W}H'_p\supset\cdots H'_2\supset H'_1\supset H'_0\simeq
H_0,
\end{equation}
where in the last step we used the usual Hilbert space isomorphism
between $H_0$ and its dual $H'_0$.

Going over to the sequence space, $H'_p$ corresponds to
\begin{equation}\label{eq:defFpprime}F_{-p}\stackrel{\rm def}{=}
\{(x_n)_{n\in W}: \sum_{n\in W}\lambda_n^{-2p}x_n^2<\infty\}
\end{equation}
The element $y\in F_{-p}$ corresponds to the linear functional on
$F_p$ given by
$$x\mapsto \sum_{n\in W}x_ny_n$$
which, by Cauchy-Schwarz, is well-defined and does define an element
of the dual $F'_p$ with norm equals to $|y|_{-p}$.

\subsection{Gaussian measure in infinite dimensions}\label{ss:gmi}

Consider now the product space ${\mathbb{R}}^W$, along with the
coordinate projection maps $${\hat X}_j: {\mathbb{R}}^W\to
{{\mathbb{R}}}:x\mapsto x_j$$ for each $j\in W$. Equip
${\mathbb{R}}^W$ with the product $\sigma$--algebra, i.e. the
smallest $\sigma$--algebra with respect to which each projection map
${\hat X}_j$ is measurable.   Kolmogorov's theorem on infinite
products of probability measures provides a
 probability measure $\nu$ on the product $\sigma$--algebra  such that
 each function ${\hat X}_j$, viewed as a random variable,
  has standard Gaussian distribution. Thus,
 $$\int_{{\mathbb{R}}^W}e^{it{\hat X}_j}\,d\nu=e^{-t^2/2}$$
 for $t\in{\mathbb{R}}$,  and every $j\in W$. The measure $\nu$ is the
product of the
 standard Gaussian measure $e^{-x^2/2}(2\pi)^{-1/2}dx$ on each
 component ${{\mathbb{R}}}$ of the product space ${\mathbb{R}}^W$.

Since, for any   $p\geq 1$, we have
$$\int_{{{\mathbb{R}}^W}}\sum_{j\in
W}\lambda_j^{-2p}x_j^2\,d\nu(x)=\sum_{j\in W}\lambda_j^{-2p}<\infty,$$
it follows that $\nu(F_{-p})=1$ for all $p\geq 1$. Thus $\nu(F')=1$.

We can, therefore, transfer  the measure $\nu$ back to $\Hstar$ ,
obtaining a probability measure $\mu$ on the $\sigma$--algebra of
subsets of $\Hstar$  generated by the maps
$$
{\hat e}_j:\Hstar\to {\mathbb{R}}:f\mapsto f(e_j),
$$
where $\{e_j\}_{j\in W}$ is the orthonormal basis of $H_0$ we
started with (note that each $e_j$ lies in $\Hnucl=\bigcap_{p\geq 0}H_p$).
This is clearly the $\sigma$--algebra generated by the weak topology
on $\Hstar$  (which happens to be equal also to the $\sigma$--algebras
generated by the strong/inductive-limit topology \cite{JB06}).

Specialized to the example $H_0=L^2({\mathbb{R}})$, and
$A=-\frac{d^2}{dx^2}+\frac{x^2}{4}+\frac{1}{2}$, we have the
standard Gaussian measure on the distribution space
$\mathcal{S}'({\mathbb{R}})$.

The above discussion gives a simple direct description of the
measure $\mu$. Its existence is also obtainable by applying the
 Minlos theorem:

\begin{thm}[Minlos theorem] \label{T:minlos}
A complex value function $\phi$ on a nuclear space $\Hnucl$  is the characteristic function of a unique probability measure $\nu$ on $\Hstar$ , i.e.,
\[
\phi(v) = \int_{E'} e^{i\ip<x,y>} \, d\nu(x) = \int_{E'} e^{i\hat{y}(x)} \, d\nu(x), \quad y\in \Hnucl
\]
if and only if 
$\phi(0) = 1 $, $\phi$ is continuous, and $\phi$ is positive definite.
\end{thm}
For a proof of the Minlos theorem refer to \cite{GV64}.
Applying the Minlos theorem to the characteristic function $\phi(y) = e^{-\frac{1}{2}\norm|y|_0^2}$ gives us the standard Gaussian measure $\mu$ we just constructed.

 There is also the useful standard
setting of Abstract Wiener Spaces for Gaussian measures introduced
by L. Gross (see the account in Kuo \cite{hk75}).

To summarize, we can state the starting point of much of
infinite-dimensional distribution theory (white noise analysis):
Given a real, separable Hilbert space $H_0$ and a positive
Hilbert-Schmidt operator $A$ on $H_0$, we have constructed a nuclear
space $\Hnucl$  and a unique probability  measure $\mu$ on the Borel
$\sigma$--algebra of the dual $\Hstar$  such that there is a linear map
$$H_0\to L^2(\Hstar,\mu):x\mapsto {\hat x},$$
satisfying
$$\int_{\Hstar}e^{it{\hat x}}\,d\mu=e^{-t^2|x|_0^2/2},$$
for every real $t$ and $x\in H_0$. This Gaussian measure $\mu$ is
often called the {\it white noise measure} and forms the background
measure for white-noise distribution theory.

\section{White Noise Distribution Theory}\label{sec:wndt}

We can now develop the ideas of the preceding section further to
construct a space of test functions   over the dual space $\Hstar$,
where $\Hnucl$  is the nuclear space related to a
 real separable Hilbert space $H_0$ as in the
  discussion in Section \ref{ss:taf}. We will use the notation,
  and in particular the spaces $H_p$, from Section \ref{ss:taf}.

The symmetric Fock space  ${\mathcal F}_s(V)$ over a Hilbert space
$V$ is the subspace of symmetric tensors in the completion of the
tensor algebra $T(V) $ under the inner--product given by
\begin{equation}\label{eq:defip}
\ip<a,b>_{T(V)} = \sum_{n = 0}^{\infty} n!
\ip<a_n,b_n>_{V^{\otimes n}},
\end{equation} where $a=\{ a_n \}_{n\geq 0}, b=
\{b_n\}_{n\geq 0}$ are elements of $T(V)$ with  $a_n, b_n$ in the
  tensor power $V^{\otimes n}$.  Then we have
 \begin{equation}\label{eq:chain}
 {\mathcal F}_s({ \Hnucl})\stackrel{\rm def}{=}
 \bigcap_{p\geq 0}{\mathcal F}_s(H_p)\subset\cdots\subset
  {\mathcal F}_s(H_2)\subset {\mathcal F}_s(H_1)\subset {\mathcal
  F}_s(H_0).
  \end{equation}
  Thus, the pair $\Hnucl \subset H_0$  give  rise to a corresponding
  pair by taking symmetric Fock spaces:
  \begin{equation}\label{eq:symFock}
 {\mathcal F}_s({ \Hnucl}) \subset {\mathcal
  F}_s(H_0).
  \end{equation}

  \subsection{\WI\ Isomorphism} \label{ss:GaussM}

  In infinite dimensions the role of Lebesgue measure is played
  by Gaussian measure $\mu$.
    There is a standard unitary isomorphism, the {\sl Wiener-It\^o
    isomorphism} or wave-particle duality map,
    which identifies the
  complexified Fock space ${\mathcal F}_s(H_0)_c$ with
  $L^2(\Hstar,\mu)$. This is uniquely
  specified by
  \begin{equation}\label{eq:WienerIto}
  I:{\mathcal F}_s(H_0)_c\to L^2(\Hstar,\mu):
  {\rm Exp}(x)\mapsto e^{{\hat
  x}-\frac{1}{2}|x|_0^2}
  \end{equation}
  where $x \in \Hnucl$  and
  $${\rm Exp}(x)=\sum_{n= 0}^{\infty}\frac{1}{n!}x^{\otimes n}.$$
   Indeed, it is readily checked that $I$
  preserves inner--products (the inner--product is as described in
  (\ref{eq:defip})). Using $I$, for each $\fockp$ with
  $p \geq 0$, we have the corresponding space $[H]_p
  \subset \LtwoH$ with the
   norm $\Vert \cdot \Vert_p$ induced by the norm on
   the space $\fockp_c$. 
   The chain of spaces   (\ref{eq:chain}) can be transferred into a chain   of function
    spaces:
\begin{equation}\label{eq:chainfunct}
 \TE=\bigcap_{p\geq 0}[H]_p\subset\cdots\subset
  [H]_2\subset [H]_1\subset [H]_0=L^2(\Hstar,\mu).
\end{equation}
     Observe that $[\Hnucl]$ is a nuclear space with topology induced by the norms $\{\Vert \cdot \Vert_p \, ; \, p = 0,1,2,\dots \}$.
  Thus, starting with the pair $ \Hnucl \subset H_0$
one obtains a corresponding pair
$[\Hnucl]\subset L^2(\Hstar,\mu).$

As before, the identification of $H'_0$ with $H_0$ leads to a complete chain
\begin{equation}\label{eq:fullchain}
 \Hnucl =\bigcap_{p\geq 0}H_p\subset\cdots \subset
 H_1\subset H_0\simeq H_{-0}\subset
 H_{-1}\subset\cdots\subset\bigcup_{p\geq 0}H_{-p}=\Hstar.
 \end{equation}
 In the same way we have a chain for the `second quantized' spaces
 ${\mathcal F}_s(H_q)_c\simeq [H]_q$. The unitary isomorphism $I$ extends to unitary
 isomorphisms
\begin{equation}\label{eq:Iextendp}
I:{\mathcal F}_s(H_{-p})_c\to [H]_{-p}\stackrel{\rm def}{=}
 [H]_{p}'\subset
\GE,
 \end{equation}
 for all $p\geq 0$. In more detail, for $a\in {\mathcal F}_s(H_{-p})_c$
 the distribution $I(a)$ is specified by
 \begin{equation}\label{eq:Ia}
 \langle I(a), \phi \rangle =\langle a,  {I^{-1}(\phi)}\rangle,
 \end{equation}
 for all $\phi \in \TE$. On the right side here we have the
 pairing of ${\mathcal F}_s(H_{-p})_c$ and ${\mathcal F}_s(H_{p})_c$
 induced by the duality pairing of $H_{-p}$ and $H_p$; in
 particular, the pairings above are complex bilinear (not
 sesquilinear).

 \subsection{Properties of test functions}\label{ss:prop}
The following theorem summarizes the properties of $\TE$ which are commonly used.
The results here are standard (see, for instance,  the monograph
\cite{HK96} by Kuo), and we compile them here for ease of
 reference.

\begin{theorem}\label{th:testfunct} Every function in  $\TE$
is $\mu$-almost-everywhere equal to a unique continuous function on
$\Hstar$. Moreover, working with these continuous versions,
\begin{enumerate}
\item $\TE$ is an algebra under pointwise
 operations; \label{I:closed}
 \item pointwise addition and multiplication are continuous
 operations $\TE \times  \TE\to \TE$; \label{I:cont}
 \item for any $\phi\in \Hstar$, the evaluation map
 $$\delta_\phi:\TE \to \R: F \mapsto F(\phi)$$
 is continuous; \label{I:delta}
 \item the exponentials $e^{i{\hat x}-\tfrac{1}{2}|x|_0^2}$, with $x$
 running over $\Hnucl$, span a dense   subspace of
 $\TE$. \label{I:dense}
 \end{enumerate}
\end{theorem}

 A complete characterization of the space $\TE$ was
obtained by Y. J. Lee (see the account in Kuo \cite[page 89]{HK96}). The test functions in $\TE$ also have a useful growth condition which can be imposed on them.

\begin{thm} \label{T:testBound}
Let $\phi \in \TE$. The $\phi$ satisfies the following growth condition for any $p \geq 0$, 
\[
|\phi(x)| \leq K_p \exp\left[\frac{1}{2} \norm|x|_{-p}^2 \right], \, \qquad x\in H_p'
\]
where $K_p$ is a constant depending on the choice of $p$.
\end{thm}
A proof of this exponential bound can be found in \cite{HK96} (Theorem 6.8 page 55). 

 \subsection{The Segal--Bargmann Transform}\label{ss:SB}
  An important tool for studying test functions and distributions in the white noise setting is the Segal--Bargmann transform. The
Segal--Bargmann transform takes a function $F \in \LtwoH$ to the
function
$SF$ on the complexified space $\Hnucl_c$ given by
\begin{equation}\label{eq:Sfzclass}
SF(z)=\int_{{ H}'}  e^{{\tilde z}-z^2/2}F\,d\mu, \qquad z \in \Hnucl_c
\end{equation}
with notation as follows: if $z=a+ib$, with $a,b\in {\mathcal H}$
then
 \begin{equation}\label{eq:defxz} {\tilde z}(x)  \stackrel{\rm
def}{=}zx\stackrel{\rm def}{=}\langle x,a\rangle+i\langle
x,b\rangle, \qquad \text{for } x \in \Hstar
\end{equation}
 and $z^2=zz$, where the product $zu$ is specified through
 \begin{equation}\label{eq:defzu}zu\stackrel{\rm def}{=}
 \langle a,s\rangle-\langle
b,t\rangle + i\left(\langle a,t\rangle+\langle
b,s\rangle\right)
\end{equation} if $z=a+ib$ and $u=s+it$, where $a,b,s,t\in \Hnucl$.

 Let $\mu_c$ be the Gaussian
measure $\Hstar_c$ specified by the requirement that
\begin{equation}\int_{\Hstar_c}e^{ax+by}\,d\mu_c(x+iy)=e^{(a^2+b^2)/4}
\end{equation}
for every $a,b\in \Hnucl$. For convenience, let us introduce the
renormalized exponential function $c_w=e^{{\tilde w}-w^2/2}\in
L^2(\Hstar,\mu)$ for all $w\in \Hnucl_c$.   It is readily checked that
for any $w \in \Hnucl_c$
\begin{equation}\label{eq:Scwz}
[Sc_w](z)=e^{wz}, \qquad \text{ for all } z \in \Hnucl_c.
\end{equation}
 Thus we may
take $Sc_w$ as a function on $\Hstar_c$ given by
$Sc_w=e^{{\tilde w}}$
where now ${\tilde w}$ is a function on $\Hstar_c$ in the
natural way.  Then $Sc_w\in L^2(\Hstar_c,\mu_c)$ and one has
\[
\left\langle Sc_w,Sc_u\right\rangle_{L^2(\mu_c)}=\langle
c_w,c_u\rangle_{L^2(\mu_c)}=e^{w{\overline u}}.
\]
 This shows that $S$
provides an isometry from the linear span of the exponentials $c_w$
in $\LtwoH$ onto the linear span of the complex exponentials
$e^{\tilde{w}}$ in $L^2(\Hstar_c,\mu_c)$. Passing to the closure one
obtains the {\bf Segal--Bargmann} unitary isomorphism
$$S:L^2(\Hnucl,\mu)\to Hol^2(\Hstar_c,\mu_c)$$
where  $Hol^2(\Hstar_c,\mu_c)$ is the closed linear span of the
complex exponential functions $e^{\tilde{w}}$ in $L^2(\Hstar_c,\mu_c)$.

An explicit expression for $SF(z)$ is suggested by
(\ref{eq:Sfzclass}). For any $\phi \in \TE$ and $z\in
\Hstar_c$, we have
\begin{equation}\label{eq:defSfzIExp}
(S\phi)(z)=\left\langle I \! \left({\rm Exp}(z)\right),\phi\right\rangle
\end{equation}
where the right side is the evaluation of the distribution
$I \! \left({\rm Exp}(z)\right)$ on the test function $\phi$. Indeed it may
be readily checked that if $S\phi(z)$ is defined in this way then
$[Sc_w](z)=e^{wz}$.

In view of (\ref{eq:defSfzIExp}), it natural to extend the Segal-Bargmann transform to distributions:  for $\Phi\in \TE'$, define $S\Phi$ to be the function on $\Hnucl_c$ given by
\begin{equation}\label{def:SPhi}
S\Phi(z)\stackrel{\rm def}{=}\left\langle  \Phi,I \!
 \left({\rm Exp}(z)\right) \right\rangle, \qquad z \in \Hnucl_c
\end{equation}

One of the many applications of the the $S$--transform includes its usefulness in  characterizing generalized functions in $\Hdist$. 

\begin{theorem} [Potthoff--Streit] \label{T:PS}
Suppose a function $F$ on $\Hnucl_c$ satisfies:
\begin{enumerate}
\item For any $z,w \in \Hnucl_c$, the function $F(\ga z+w)$ is an entire function of $\ga \in \C$. \label{L:one}
 \item There exists nonnegative constants $A, p,$ and $C$ such that \label{L:two}
\[
\vert F(z) \vert \leq C e^{A\norm|z|_p^2} \qquad \text{ for all } z \in \Hnucl_c.
\]
\end{enumerate}
Then there is a unique generalized function $\Phi \in \Hdist$ such that $F = S\Phi$. Conversely, given such a $\Phi \in \Hdist$, then $S\Phi$ satisfies \emph{(\ref{L:one})} and \emph{(\ref{L:two})} above.
\end{theorem}

For a proof see Theorem 8.2 in Kuo's book \cite{HK96} on page 79.

The $S$-transform can also aid us in determining convergence in $\GE$.
\begin{thm} \label{T:conv}
Let $\Phi_n \in \GE$ and $F_n = S\Phi_n$. Then $\Phi_n$ converges strongly in $\GE$ if and only if the following conditions are satisfied:
\begin{enumerate}
\item $\lim_{n\to\infty} F_n(z)$ exists for all $z \in \Hnucl_c$.
 \item There exists nonnegative constants $A, p,$ and $C$ such that 
\[
\vert F_n(z) \vert \leq C e^{A\norm|z|_p^2}, \text{ for all } n \in N, z \in \Hnucl_c.
\]
\end{enumerate}
\end{thm}

For a proof see Kuo's book \cite{HK96} (Page 86, Theorem 8.6).

\subsection{Translation of the Gaussian Measure} The Gaussian measure $\mu$ on $\Hstar$ and its translation $\mu(\cdot - \xi)$ are related via the \ST\ when $\xi \in \Hnucl$ \cite{NO94}. Observe the following:
\begin{prop} \label{P:translation}
The Gaussian Measure $\mu$ is quasi-invariant under the translation by any $\xi \in \Hnucl$ and the Radon-Nikodym derivative is given by
\[
\frac{d\mu(\cdot -\xi)}{d\mu} = e^{\ip<\cdot,\xi> - \frac{1}{2}\ip<\xi,\xi>}.
\]
\end{prop}

\begin{proof}
Suppose $x \in \Hnucl$ and  consider the measure given by 
\[
\gl(A) = \int_{A} e^{\ip<x,\xi> - \frac{1}{2}\ip<\xi,\xi>/2} \, d\mu(x).
\]
We compute the characteristic equation
\[
\hat{\gl}(y) = \int_{\Hstar} e^{i\ip<x,y>} e^{\ip<x,\xi> - \frac{1}{2}\ip<\xi,\xi>} \, d\mu(x).
\]
Then using the characteristic equation for the Gaussian measure we have that the above gives us
\[
\hat{\gl}(y) = e^{\frac{1}{2}\ip<\xi +iy,\xi+iy> - \frac{1}{2}\ip<\xi,\xi>} = e^{i\ip<\xi,y> -\ip<y,y>/2}.
\]
On the other hand, since
\[
 \int_{\Hstar} e^{i\ip<x,y>} \, d\mu(x-\xi) = \int_{\Hstar} e^{i\ip<x+\xi, y>} \, d\mu(x)=e^{i\ip<\xi,y>-\ip<y,y>/2}
 \]
 we have that
 \[
 \hat{\gl}(y) = \int_{\Hstar} e^{i\ip<x,y>} \, d\mu(x-\xi)
 \]
 and hence $
\frac{d\mu(\cdot -\xi)}{d\mu} = e^{\ip<\cdot,\xi> - \frac{1}{2}\ip<\xi,\xi>}
$.
\end{proof}

\subsection{Translation Operator} An important operator acting on the space of test functions is the translation operator $T_y$ with $y \in \Hstar$. 

\begin{Def}
For any $y \in \Hstar$ the translation operator $T_y$ on $\TE$ is defined by
\[
T_y\phi(x) = \phi(x+y).
\]
\end{Def}

Since the Gaussian measure is not translation invariant this operator is more intricate than it first appears. 
The properties of this operator are summarized in the following theorem.
\begin{thm} \label{T:trans}
For any $y \in \Hstar$, the translation operator $T_y$ is continuous from $\TE$ into itself. Moreover, if $y \in H_p'$ and $q >p$ satisfies $\gl_1^{2(q-p)} > 2$ then for all $\phi \in \TE$,
\[
\Norm|T_y\phi|_p \leq \Norm|\phi|_q(1-\tfrac{2}{\gl_1^{2(q-p)}})\exp \left[ \tfrac{1}{\gl_1^{2(q-p)}}\norm|y|_{-p}^2 \right].
\]
\end{thm}

A proof of this can be found in the book by Kuo \cite{HK96} (page 138, Theorem 10.21).

\section{Gaussian Measure on an Affine Subspace}

For a subspace $W$ of $\Rn$ and $a \in W^{\perp}$ we have the Gaussian measure on $a+W$ given by:
\[
\int_{a+W} e^{i\ip<x,y>} \, d\mu_{a+W}(x) = \int_{a+W} e^{i\ip<x,y>} e^{-\fr(1/2)\norm|x-a|^2} \, \fr(dx/{(2\pi)}^{\dim W/2}) = e^{i\ip<a,y> - \fr(1/2)\ip<y_W,y_W>}
\] 
where $y \in \Rn$ and $y_W$ is the projection of $y$ onto $W$. We now describe how such a measure can be constructed in white noise setting.
Of course, the Gaussian measure cannot live on $H_0$ or $a+W$. However, just as we used the Minlos theorem to form the Gaussian measure $\mu$ on $\Hstar$  (which we think of as the Gaussian measure on $H_0$), we can again use the Minlos theorem to form the Gaussian measure for the affine subspace $a+W$.

\subsection{Gaussian Measure on $a+V$} For a vector $a \in H_0$ and a subspace $V$ of $H_0$ we can use the Minlos theorem to find that there is a measure $\mu_{a+V}$ on $\Hstar$  with 
\begin{equation} \label{eq:uaVchar}
\int_{\Hstar} e^{i\ip<x,y>} \, d\mu_{a+V}(x) = e^{i\ip<a,y> - \fr(1/2)\ip<y_V, y_V>} 
\end{equation}
for any $y \in \Hnucl$. This measure $\mu_{a+V}$ is the \emph{Gaussian measure for the affine subspace $a +V$}. 
This measure was originally constructed in \cite{jB07}.

\subsection{Hida Measure} The Gaussian measure $\mu_{a+V}$ is a special type of measure known as a Hida measure. In this section we define the notion of Hida measure and give an overview of some its properties. 

\begin{Def}
A measure $\gv$ on $\Hstar$  is called a \emph{Hida measure} if $\phi \in L^1(\gv)$ for all $\phi \in \TE$ and the linear functional
$$ \phi \mapsto \int_{\Hstar} \phi(x) \, d\gv(x)$$
is continuous on $\TE$.
\end{Def}

We say that a generalized function $\Phi \in \GE$ is \emph{induced} by a Hida measure $\gv$ if for any $\phi \in \TE$ we have
\[
\pair<\Phi, \phi> = \int_{\Hstar} \phi(x) \, d\gv(x).
\]
The following theorem characterizes those generalized functions which are induced by a Hida measure.
\begin{theorem} \label{T:Hida}
Let $\Phi \in \GE$. Then the following are equivalent:
\begin{enumerate}
\item For any nonnegative $\phi \in \TE$, $\pair<\Phi, \phi> \geq 0$. 
\item The function $\T (\Phi)(x) = \pair<\Phi, e^{i \ip<\cdot,x>}>$ is positive definite on $\Hnucl$.
\item $\Phi$ is induced by a Hida measure.  
\end{enumerate}
\end{theorem}

A proof of this theorem can be found in \cite{HK96} (page 320, Theorem 15.3).

\begin{cor} \label{C:Hida}
Let $\gv$ be a finite measure on $\Hstar$  such that for any $x \in \Hnucl$
\[
\pair<\Phi, e^{i \ip<\cdot,x>}> = \int_{\Hstar} e^{i \ip<y,x>} \, d\gv(y)
\]
for some $\Phi \in \GE$. Then $\Phi$ is induced by $\gv$.
\end{cor}
\begin{proof}
Since $\pair<\Phi, e^{i \ip<\cdot,x>}> = \int_{\Hstar} e^{i \ip<y,x>} \, d\gv(y)$ it is clear that $\pair<\Phi, e^{i \ip<\cdot,x>}>$ is positive definite. So we can apply \thmref{T:Hida} to get a finite measure $m$ which is induced by $\Phi$. Hence for all $\phi \in \TE$,
\[
\pair<\Phi, \phi> = \int_{\Hstar} \phi \, dm.
\]
Letting $\phi = e^{i \ip<\cdot,x>}$ in the above equation, we see that the characteristic functions for $m$ and $\gv$ are identical. Therefore $m = \gv$ and we have that $\Phi$ is induced by $\gv$.
\end{proof}

Here is another useful theorem which characterizes Hida measures.
\begin{thm} \label{T:Hidabnd}
A measure $\nu$ on $\Hstar$ is a Hida measure if and only if $\nu$ is supported in $H_p'$ for some $p \geq 1$ and
\[
\int_{H_p'} \exp\left[\tfrac{1}{2}|x|_{-p}^2\right] \, d\nu(x) < \infty.
\]
\end{thm}

For a proof of this refer to Kuo's book \cite{HK96} (page 333, Theorem 15.17).

\subsection{Definition of the distribution $\dV[a+V]$} We now prove that $\mu_{a+V}$ is a Hida measure and develop the corresponding distribution $\dV[a+V]$ which we think of as the delta function for the affine subspace $a+V$ \cite{jB07}.
Observe the effect of $\mu_{a+V}$ on the renormalized exponential $e^{\ip<\cdot, z> - \fr(1/2) \ip<z,z>}$,
\begin{align*}
\int_{\Hstar} e^{\ip<x,z> - \fr(1/2) \ip<z,z>} \, d\mu_{a+V}(x)  
&= e^{-\ip<z,z>} \int_{\Hstar} e^{\ip<x,z>} \, d\mu_{a+V}(x)  \\
&= e^{-\ip<z,z>} e^{\ip<a,z> + \fr(1/2) \ip<z_V,z_V>}  \\
&=e^{\ip<a,z> - \fr(1/2) \ip<z_{\Vperp},z_{\Vperp}>}.
\end{align*}
Although $\dV[a+V]$ was originally developed for $a \in H_0$ we could also take $a \in H_p'$.
Let the function $F(z)$ denote the result from the calculations above. That is,
\begin{equation} \label{E:F}
F(z) = e^{\ip<a,z>- \fr(1/2) \ip<z_{\Vperp},z_{\Vperp}>}
\end{equation} 
We will show that $F(z)$ satisfies properties (1) and (2) of \thmref{T:PS}. 

For property (1) consider $F(\ga z + w)$ where $z,w \in \Hnucl_c$ and $\alpha \in \C$. Then notice that
\begin{align*}
F(&\ga z + w) =  e^{\ip<a, \ga z + w> - \tfr(1/2) \ip<\ga z_{\Vperp} + w_{\Vperp},\ga z_{\Vperp} + w_{\Vperp}> } \\
&= \exp[\ga \ip<a,z> + \ip<a,w> - \tfr(1/2) ( \ga^2\ip<z_{\Vperp}, z_{\Vperp}> + 2\ga\ip<z_{\Vperp}, w_{\Vperp}>+\ip<w_{\Vperp}, w_{\Vperp}>)] \\
&= e^{-\tfr(\ga^2/2) \ip<z_{\Vperp} , z_{\Vperp}>} e^{ \ga(\ip<a,z> -\ip<z_{\Vperp}, w_{\Vperp}>)} e^{\ip<a,w> - \tfr(1/2) \ip<w_{\Vperp}, w_{\Vperp}>}
\end{align*}
which is an  entire function of $\ga \in \C$.

Now for property (2) of \thmref{T:PS} we write $z$ as $z = x+iy$ with $x,y \in \Hnucl$ and observe that
\begin{align*}
\vert F(z) \vert  &= | e^{\ip<a,z> - \fr(1/2)\ipc<z_{\Vperp},z_{\Vperp}> }| \\
&= | e^{\ip<a,x+iy> - \fr(1/2)\ipc<x_{\Vperp}+iy_{\Vperp},x_{\Vperp}+iy_{\Vperp}> } | \\
&= e^{\ip<a,x>} e^{-\fr(1/2)\norm|x_{\Vperp}|_0^2+\fr(1/2)\norm|y_{\Vperp}|_0^2} \\
&\leq e^{\ip<a,x>}  e^{\fr(1/2)\norm|z_{\Vperp}|_0^2} \\
&\leq e^{\norm|a|_{-p}\norm|x|_{p}}e^{\fr(1/2)\norm|z|_0^2} \\
&\leq e^{\fr(1/2) \norm|a|_{-p}^2 + \fr(1/2) \norm|z|_p^2} e^{\fr(1/2)\norm|z|_p^2} \qquad \text{by Young's Inequality} \\
&\leq e^{\fr(1/2) \norm|a|_{-p}^2} e^{\fr(3/2)\norm|z|_p^2}. 
\end{align*}
So property (2) of \thmref{T:PS} is satisfied.

Therefore by \thmref{T:PS} there exist some $\Phi \in \Hdist$ such that $S(\Phi)(z) = F(z)$.
Then by \corref{C:Hida} we have that  for $a \in H_0$, $\Phi$ is induced by $\mu_{a+V}$. We simply denote this $\Phi$ by $\dV[a+V]$. This leads us to the following definition: \cite{jB07}

\begin{Def}  \label{D:dV}
The \emph{delta function for the affine subspace} $a+V$ is the distribution in $\GE$ induced by the Hida measure $\mu_{a+V}$. We denote this generalized function by $\dV[a+V]$.
\end{Def} 

Thus for any test function $\phi \in \TE$ we have 
\[
\pair<\dV[a+V], \phi> = \int_{\Hstar} \phi \, d\mu_{a+V}.
\]

\subsection{\ST of \dV[a+V]}
Using the definition of the distribution $\dV[a+V]$ we can directly compute its \ST. By the calculations directly preceding \eqref{E:F} we have 
\begin{equation} \label{eq:Strans}
S(\dV[a+V])(z)  = \SdV \qquad \text{for} \quad z \in \Hnucl_c.
\end{equation}
Using this framework of the Hida measure $\mu_{a+V}$ and the corresponding distribution $\dV[a+V]$ we have the following intuitive theorem:
\begin{theorem}\label{T:shift}
Let $V$ be a subspace of $H_0$ and $a \in \Vperp$, then for any $\phi \in \TE$ we have
\begin{equation}\label{eq:shift}
\int_{\Hstar} \phi(x) \, d\mu_{a+V}(x) = \int_{\Hstar} \phi(x+a) \, d\mu_{V}(x)
\end{equation}
\end{theorem}

\begin{proof}
First we take the special case where $\phi(x) = e^{i\ip<x,\xi>}$ for some $\xi \in \Hnucl$. Then we have for the left hand side
\[
\int_{\Hstar} \phi(x) \, d\mu_{a+V}(x) =\int_{\Hstar} e^{i\ip<x,\xi>} \, d\mu_{V}(x) = e^{i\ip<a,\xi> -\frac{1}{2}\ip<\xi_V, \xi_V>}
\]
and for the right hand side
\begin{align*}
\int_{\Hstar} \phi(x+a) \, d\mu_{V}(x) &=\int_{\Hstar} e^{i\ip<x+a,\xi>} \, d\mu_{V}(x) \\
& = e^{i\ip<a,\xi>}\int_{\Hstar} e^{i\ip<x,\xi>} \, d\mu_{V}(x) = e^{i\ip<a,\xi> -\frac{1}{2}\ip<\xi_V, \xi_V>}
\end{align*}
Thus we have that \eqref{eq:shift} agrees on the linear span of $\{ e^{i\ip<x,\xi>} \, ; \, \xi \in \Hstar \}$. 

For any arbitrary $\phi \in \TE$ take a sequence $\phi_n$ in the linear space of $\{ e^{i\ip<x,\xi>} \, ; \, \xi \in \Hstar \}$ such that $\phi_n$ converges to $\phi$ in $\TE$. Then we have
\begin{align*}
\int_{\Hstar} \phi(x) \, d\mu_{a+V}(x) &= \pair<\phi, \dV[a+V]> \\
&= \lim_{n\to\infty} \pair<\phi_n, \dV[a+V]> \\
&= \lim_{n\to\infty} \int_{\Hstar} \phi_n(x) \, d\mu_{a+V}(x)\\
&= \lim_{n\to\infty} \int_{\Hstar} \phi_n(x+a) \, d\mu_{V}(x) \\
&=\lim_{n\to\infty} \pair<T_a(\phi_n), \dV[V]> \\
&= \pair<T_a(\phi), \dV[V]> \text{ using the continuity of $T_a$} \\
&= \int_{\Hstar} \phi(x+a) \, d\mu_{V}(x)
\end{align*}
giving us the desired result.
\end{proof}
Now we prove a convenient and somewhat expected property of convergence amongst these delta functions on an affine subspace.

\begin{prop} \label{P:deltaConv}
Let $\{x_n\}$ be a sequence in $H_p'$ converging to $x$  and suppose $\{S_n\}$ is a sequence of subspaces of $H_0$ converging to a subspace $S$, in the sense that for any $v \in H_0$, we have $v_{S_n}$ converges to $v_{S}$ in $H_0$.
Then the generalized functions $\dV[x_n + S_n]$ converges strongly to  $\dV[x+S]$ in $\GE$.
\end{prop}

\begin{proof}
We will apply \thmref{T:conv}. To see that the conditions of \thmref{T:conv} are satisfied notice that for $z \in \Hstar_c$ we have
\begin{align*}
\lim_{n\to\infty} S(\dV[x_n + S_n])(z) &= \lim_{n\to\infty}  \pair<\dV[x_n + S_n], e^{\ip<\cdot,z>-\ip<z,z>}>  \\
&= \lim_{n\to\infty} e^{\ip<x_n, z> - \frac{1}{2}\ip<z_{S_n^{\perp}}, z_{S_n^{\perp}}>} \text{ by \eqref{eq:Strans} }\\ 
&= e^{\ip<x, z> - \frac{1}{2}\ip<z_{S^\perp}, z_{S^\perp}>} \\
&= \pair<\dV[x+S], e^{\ip<\cdot,z>-\frac{1}{2}\ip<z,z>}>\\
&= S(\dV[x+S])(z).
\end{align*}
For the second condition of \thmref{T:conv} notice that
\begin{align*}
S(\dV[x_n + S_n])(z) &= e^{\ip<x_n, z> - \frac{1}{2}\ip<z_{S_n^{\perp}}, z_{S_n^{\perp}}>} \\
                     &\leq e^{\norm|x_n|_{-p} \norm|z|_p}e^{\frac{1}{2}\norm|z_{S_n^{\perp}}|_0^2 }\\
                     &\leq e^{\norm|x|_{-p} \norm|z|_p}e^{\frac{1}{2}\norm|z|_0^2 }\\
                     &\leq e^{\frac{1}{2}\norm|x|_{-p}^2 + \frac{1}{2}\norm|z|_p^2}e^{\frac{1}{2}\norm|z|_0^2 } \text{ by Young's Inequality} \\
                     &\leq  e^{\frac{1}{2}\norm|x|_{-p}^2 + \frac{1}{2}\norm|z|_p^2}e^{\frac{1}{2}\norm|z|_p^2 }\\
                     &= e^{\frac{1}{2}\norm|x|_{-p}^2 }e^{\frac{3}{2}\norm|z|_p^2 }.
\end{align*}
\end{proof}

\section{Gauss Radon Transform in Infinite Dimensions} 
We begin by constructing the Radon--Gauss Transform in $\R^n$. Recall that a hyperplane in $\R^n$ can be represented using a unit normal vector $v \in \R^n$ and a number $\ga \in \R$ by way of
\[
\ga v + v^{\perp}.
\]
The  probability density function for the standard Gaussian measure $\mu_{\ga v+ v^{\perp}}$ on $\ga v+ v^{\perp}$ is given by
\[
d\mu_{\ga v + v^{\perp}}(x) = \frac{1}{(2\pi)^{\frac{n-1}{2}}} e^{-|x-\ga v|^2/2 }\, dx
\]
where  $x \in \R^n$, but $dx$ denotes the Lebesgue measure on $\ga v + v^{\perp}$. The characteristic function of this measure is given by
\begin{equation} \label{eq:charRn}
\hat{\mu}_{\ga v + v^{\perp}}(k) = e^{i\ga\ip<k,v>-\frac{1}{2}\ip<k_{v^{\perp}}, k_{v^{\perp}}>},
\end{equation}
where $k_{v^{\perp}}$ is the orthogonal projection of $k$ onto $v^{\perp}$.

Using the measure $\mu_{a+V}$ we can construct the Gauss--Radon transform in the white noise framework. (Note that the Gauss--Radon transform was originally constructed for a similar setting in \cite{vM00}.)
\subsection{Hyperplanes in $H_0$}
In infinite dimensions we define a hyperplane as follows:
\begin{Def}
A \emph{hyperplane}  of a infinite dimensional Hilbert space $H_0$ is given by the set
\[
\ga v+ v^{\perp} = \{ \ga v + x \, ; \, x \in H_0, \ip<x,v>_0=0 \}
\]
where $\ga$ is a real number and $v$ is a non-zero unit vector in $H_0$.
\end{Def}
 For such an affine subspace the measure $\mu_{\ga v + v^{\perp}}$ has the following characteristic equation and \ST:
\begin{equation} \label{eq:charuv}
\int_{\Hstar} e^{i\ip<x,y>} \, d\mu_{\ga v + v^{\perp}}(x) = e^{i\ga\ip<v,y> - \fr(1/2)\ip<y_{v^{\perp}}, y_{v^{\perp}}>},  \quad y \in \Hnucl
\end{equation}
and
\begin{equation}
\int_{\Hstar} e^{\ip<x,z>-\ip<z,z>} \, d\mu_{\ga v + v^{\perp}}(x) = e^{\ga\ip<v,z> - \fr(1/2) \ip<z,v>^2},  \quad z \in \Hnucl_c.
\end{equation}
Notice that the above is analogous to what we have observed in $\Rn$. Using this measure $\mu_{\ga v + v^{\perp}}$ 
we can now define the Gauss--Radon transform in the white noise framework.
\begin{Def} \label{D:GR}
For a test function $\phi \in \TE$ we define the Gauss--Radon transform to be the function on the hyperplanes of $H_0$ given by
\[
G_{\phi}(\ga v + v^{\perp}) = \int_{\Hstar} \phi(x) \, d\mu_{\ga v + v^\perp}(x).
\]
\end{Def}

In \cite{vM00} Mihai and Sengupta also demonstrated that this measure can be constructed using the Kolmogorov theorem and Gaussian measures $\mu_n$ on $\Rn$  specified by
\[
\hat{\mu}_n(k) = e^{i\ga\ip<k, v_n> - \frac{1}{2}(|k|^2 - |\ip<k,v_n>|^2)}
\]
where $v_n = (\ip<v,e_1>, \dots, \ip<v,e_n>)$. Note that if  $|v_n| = 1$, then the above is the Gaussian measure  on the hyperplane $\{x \in \Rn \, ; \, \ip<v_n, x> = \ga \} = \ga v+ v^{\perp}$.

Putting these ideas together we have the following theorem
\begin{prop} \label{P:Rn}
Let $v \in \text{span}\{e_1, \dots, e_n\} \subset H_0$ and $v_n = (\ip<v,e_1>, \dots, \ip<v,e_n>) \in \Rn$. Then for any $\phi$ of the form $F(\ip<\cdot, e_1>, \dots, \ip<\cdot, e_n>)$ where $F$ is a integrable function with respect to the measure $\mu_{\ga v_n + v_n^{\perp}}$ on $\Rn$ we have
\[
G_{\phi}(\ga v + v^{\perp}) = \int_{\Hstar} \phi \, d\mu_{\ga v + v^{\perp}} = \int_{\ga v_n + v_n^{\perp}} F \mu_{\ga v_n + v_n^{\perp}}.
\]
\end{prop}

\subsection{Disintegration} Here we demonstrate a Fubini like theorem for our Gaussian measure on the affine subspace $a+V$. The theorem allows us to break up the integral into integrals over subspaces making up $V$. This will be most useful when $a+V$ is a hyperplane as in the Gauss--Radon Transform.
\begin{thm} \label{T:dis}
Let $\phi$ be a test function in $\Htest$ and consider the affine subspace $a+V$ in $H_0$. If
\[
V = S \oplus S^{\perp}
\]
where $S$ is a subspace of $H_0$, then
\begin{equation} \label{eq:Tdis}
\int_{\Hstar} \phi \, d\mu_{a+V} = \int_{\Hstar} \int_{\Hstar} \phi(x+y) \, d\mu_{a+S}(x) \, d\mu_{S^{\perp}}(y).
\end{equation}
\end{thm}

\begin{proof}
We first show that the above holds for $\phi(x) = e^{i\ip<x,\xi>}$ where $\xi \in \Hnucl$. The lefthand side of \eqref{eq:Tdis}
is simply the characteristic equation of  $\mu_{a+V}$  given by \eqref{eq:uaVchar} 
\begin{equation} \label{eq:compareAbove}
e^{i\ip<a,\xi> - \frac{1}{2}\ip<\xi_V,\xi_V>}
\end{equation}
Now for the righthand side we have
\begin{align*}
 \int_{\Hstar} \int_{\Hstar} e^{i\ip<x+y,\xi>} \, d\mu_{a +S}(x) d\mu_{S^{\perp}}(y) &= \int_{\Hstar} e^{i\ip<x,\xi>} \, d\mu_{a +S}(x) \int_{\Hstar} e^{i\ip<y,\xi>}  d\mu_{S^{\perp}}(y) \\
 &=e^{i\ip<a,\xi> - \frac{1}{2}\ip<\xi_{{S}},\xi_{{S}}>}e^{ -\frac{1}{2}\ip<\xi_{S^{\perp}},\xi_{S^{\perp}}> }\\
 &=e^{i\ip<a,\xi> - \frac{1}{2}\ip<\xi_{V},\xi_{V}>} \text{ because } V = S \oplus S^{\perp}
\end{align*}
So the above holds on the dense space given by the linear span of $\{e^{i\ip<\cdot,\xi>} \, ; \, \xi \in \Hnucl \}$.

Now for an arbitrary $\phi$, let $\phi_n$ be in the linear span of $\{e^{i\ip<\cdot,\xi>} \, ; \, \xi \in \Hnucl \}$. For the lefthand side we have
\[
\int_{\Hstar} \phi \, d\mu_{a+V} = \pair<\phi, \dV[a+V]> = \lim_{n\to \infty} \pair<\phi_n, \dV[a+V]> = \lim_{n\to \infty}  \int_{\Hstar} \phi_n \, d\mu_{a+V}
\]
using the relationship between the measure  $\mu_{a+V}$ and the distribution $\dV[a+V]$. The last term in the above equality is equal to
\begin{equation} \label{eq:disAbove}
 \lim_{n\to \infty} \int_{\Hstar} \int_{\Hstar}\phi_n(x+y) \, d\mu_{a +S}(x) d\mu_{S^{\perp}}(y)
\end{equation}
If we can pass the limit inside the integral then the proof will be complete. We will work inside out. First note that since 
$\mu_{S^{\perp}}$ is a Hida measure  by \thmref{T:Hidabnd} for some  $p \geq 1$ we have that  $\mu_{S^{\perp}}(H_{p}') = 1$ and 
\begin{equation} \label{eq:aboveThree}
\int_{H_p'} \exp\left[\tfrac{1}{2}|y|_{-p}^2\right] \, d\mu_{S^{\perp}}(y) < \infty 
\end{equation}
Thus we can rewrite the righthand side of \eqref{eq:Tdis} as
\begin{equation} \label{eq:disAboveTwo}
 \int_{H_p'} \int_{\Hstar}\phi(x+y) \, d\mu_{a +S}(x) d\mu_{S^{\perp}}(y).
\end{equation}
Working inside out the inside part of the above integral can be written
\[
\int_{\Hstar}\phi(x+y) \, d\mu_{a +S}(x) = \pair<T_y\phi, \dV[a +S]>.
\]
with $y \in H_p'$. Since $\dV[a +S]$ is in $\GE$ and $T_y$ is continuous from $\TE$ into itself we have that
\begin{align*}
\int_{\Hstar}\phi(x+y) \, d\mu_{a +S}(x) &= \pair<T_y\phi, \dV[a +S]> \\
&= \lim_{n\to\infty} \pair<T_y\phi_n, \dV[a +S]> = \lim_{n\to\infty} \int_{\Hstar}\phi_n(x+y) \, d\mu_{a +S}(x)
\end{align*}
Thus \eqref{eq:disAboveTwo} becomes 
\begin{align*}
\int_{H_p'} \int_{\Hstar}\phi(x+y) &\, d\mu_{a +S}(x) d\mu_{S^{\perp}}(y) \\
&= \int_{H_p'} \lim_{n\to\infty} \int_{\Hstar}\phi_n(x+y) \, d\mu_{a +S}(x) \,  d\mu_{S^{\perp}}(y).
\end{align*}
We would like to use the dominated convergence theorem to pull the limit out once more. To do this notice that
$\int_{\Hstar}\phi_n(x+y) \, d\mu_{a +S}(x)$ is measurable; here note that
\[ \int_{\Hstar}e^{i\ip<x+y,\xi> }\, d\mu_{a +S}(x) = e^{i\ip<y,\xi>}e^{i\ip<a,\xi>-\frac{1}{2}\ip<\xi_{S},\xi_{S}>} \]
which is measurable. Thus 
\[
\int_{\Hstar}\phi(x+y) \, d\mu_{a +S}(x) = \lim_{n\to\infty} \int_{\Hstar}\phi_n(x+y) \, d\mu_{a +S}(x)
\]
 is measurable. Also observe that choosing a $k$ such that $k > p$ and $\dV[a + S] \in [H_k]'$ we have
\begin{align*}
\left| \int_{\Hstar}\phi_n(x+y) \, d\mu_{a +S}(x) \right| &= \left| \pair<T_y\phi_n, \dV[a +S]> \right|\\
&\leq \Norm|T_y\phi_n|_k \Norm|\dV[a + S]|_{-k} \\
&\leq \Norm|\phi_n|_q \Norm|\dV[a + S]|_{-k} \exp\left[\tfrac{1}{2}|x|_{-p}^2\right] 
\end{align*}
using \thmref{T:trans} where $q$ is chosen to ensure that $\gl_1^{2(q-k)} > 2$. Now in the above we have  $\Norm|\phi_n|_q$ is bounded because $\phi_n \to \phi$ in $\TE$. Putting this altogether we have for some number $M$
\[
\left| \int_{\Hstar}\phi_n(x+y) \, d\mu_{a +S}(x) \right| \leq M \exp\left[\tfrac{1}{2}|x|_{-p}^2\right]
\]
and integral of the righthand side of the above using the measure $\mu_{S^{\perp}}$ is finite by \eqref{eq:aboveThree}. Therefore the dominated convergence theorem applies.
\end{proof}

For this work the above theorem proves most useful when the affine subspace is actually a hyperplane $\ga v + v^{\perp}$. Then the above gives us a means by which to decompose the Gauss--Radon transform. 

\begin{cor} \label{C:dis}
Let $\phi$ be a test function in $\Htest$ and consider the hyperplane $\ga v + v^{\perp}$ in $H_0$. If
\[
v^{\perp} = S \oplus S^{\perp}
\]
where $S$ is a subspace of $H_0$, then
\begin{equation} \label{eq:dis}
G_{\phi}(\ga v + v^{\perp}) = \int_{\Hstar} \phi \, d\mu_{\ga v + v^{\perp}} = \int_{\Hstar} \int_{\Hstar} \phi(x+y) \, d\mu_{\ga v +S}(x) \, d\mu_{S^{\perp}}(y).
\end{equation}
\end{cor}

\begin{cor}
Let $V$ be  a subspace of $H_0$. Then for any test function $\phi$ we have
\[
\int_{\Hstar} \phi(x) \, d\mu(x) = \int_{\Hstar}\int_{\Hstar} \phi(x+y) \, d\mu_V(x) \, d\mu_{V^{\perp}}(y).
\]
\end{cor}

\subsection{Coordinates} Our goal here is to show that $\int_{\Hstar} \phi(x) \, d\mu_{a+V}(x)$ essentially only depends on the ``projections'' of the $x$-values to the subspace $V$. We first need the following lemma concerning our most popular dense set.

\begin{lemma}
The linear span of $\{e^{i\ip<\cdot,\xi	>} \, ; \,  \xi \in \Hnucl \}$ is dense in $L^1(\mu_{a + V})$.
\end{lemma}

\begin{proof}
A result in \cite{vM00} (Proposition 3.4) states that the linear span of $\{e^{i\ip<\cdot,\xi	>} \, ; \,  \xi \in \Hnucl \}$ is dense in $L^2(\mu_{a + V})$. Now we simply show that $L^1(\mu_{a + V})$ is dense in $L^2(\mu_{a + V})$. Let $f \in L^1(\mu_{a + V})$ with $f$ orthogonal to $L^2(\mu_{a +V})$. Our objective is to show that $f = 0$. 

Note that for any measurable set $A$ we have that $1_A$ in in $L^1(\mu_{a + V})$ and $L^2(\mu_{a + V})$. Since $f$ is orthogonal to $L^2(\mu_{a + V})$ we must have
\[
\int_{\Hstar} 1_Af \, d\mu_{a +V} = 0 
\]
In particular for the set $\{ f \geq 0 \}$ we have that
\[
0=\int_{\{ f \geq 0 \}} f \, d\mu_{a +S} = \int_{\{ f \geq 0 \}} f^+ \, d\mu_{a +S}.
\]
Thus $f^+= 0$ almost everywhere. Similarly, we can get $f^- = 0$ almost everywhere.
\end{proof}

\begin{thm} \label{T:coor}
Let $a \in \text{span}\{e_1, \dots, e_n\} \subset H_0$ and $S$ be a subspace of $H_0$ with $S \subset \text{span}\{e_1, \dots, e_n\}$. Then if $f \in L^1(\mu_{a+S})$, we have
\[
\int_{\Hstar} f(x) \, d\mu_{a+S}(x) = \int_{\text{span}\{e_1, \dots, e_n\}} f(\ip<x,e_1>e_1+\cdots+\ip<x,e_n>e_n) \, d\mu_{a+S}(x)
\]
\end{thm}

\begin{proof}
Let $P_S$ be the projection onto the subspace $S$. Observe that for any $k > n$ we have that
\begin{align*}
\int_{\Hstar} e^{it\hat{e}_k} \, d\mu_{a + S} &= e^{i\ip<a, te_k> - \frac{1}{2} \ip<tP_Se_k, tP_Se_k> }\\
&= e^0 \\
&= \int_{\R} e^{its} \, d\gd_0(s)
\end{align*}
where $\gd_0$ is the delta measure with $\gd_0({0}) = 1$. Since the characteristic function of a random variable uniquely specifies the distribution, it follows that the random variable $\hat{e}_k$ has a distribution $\gd_0$, i.e. $\hat{e}_k$ has the constant value $0$ almost everywhere. Thus the measure of the set $\hat{e}_k^{-1}(0) = \{x \in \Hstar \, ; \, \ip<x,e_k> = 0 \}$ has full measure with respect to $\mu_{a+S}$. Therefore the set $\{\hat{e}_k \neq 0\} = \{x \in \Hstar \, ; \, \ip<x,e_k> \neq 0 \} $ has $\mu_{a+S}$--measure $0$. Hence the set
\[
\bigcup_{k = n+1}^{\infty} \{\hat{e}_k \neq 0 \}
\]
has $\mu_{a+S}$ measure $0$. Likewise the complement 
\[
\left( \bigcup_{k = n+1}^{\infty} \{\hat{e}_k \neq 0 \} \right)^c = \bigcap_{k = n+1}^{\infty} \{\hat{e}_k \neq 0 \}^c =\bigcap_{k = n+1}^{\infty} \{\hat{e}_k = 0 \} = \text{span}\{e_1, \dots, e_n\}
\]
has $\mu_{a+S}$--measure $1$. Therefore for any $f \in L^{1}(\mu_{a+S})$ we have
\begin{align*}
\int_{\Hstar} f(x) \, d\mu_{a+S}(x) &= \int_{\text{span}\{e_1, \dots, e_n\}} f(x) \, d\mu_{a+S}(x) \\
&=  \int_{\text{span}\{e_1, \dots, e_n\}} f(\ip<x,e_1>e_1+\cdots+\ip<x,e_n>e_n) \, d\mu_{a+S}(x)
\end{align*}
since $x = \ip<x,e_1>e_1+\cdots+\ip<x,e_n>e_n$ when $x \in \text{span}\{e_1, \dots e_n\}$.
\end{proof}

\section{Support Theorem for Gauss--Radon Transform}

Having the Gauss--Radon transform fully developed we take on the task of developing the Support Theorem in this setting. The Support Theorem in $\Rn$ requires that the Radon transform be zero outside of some convex compact set. The typical example in $\Rn$ are the closed balls. At some point we will to appeal the Support Theorem in $\Rn$. So the sets we consider in infinite dimensions must have ``projections'' which are convex and compact. The convexity issue is easily addressed. To have the property of compactness we desire  leads us to the following definition:
\begin{Def}
A subset $C$ of $\Hstar$ is \emph{projectively compact} if the set 
\[
C_n = \{ \vec{x}_n = (\ip<x,e_1>, \dots, \ip<x,e_n>) \, ; \, x \in C \}
\]
is compact and $\vec{x}_n \in C_n$ for all $n$ implies $x \in C$.
\end{Def}

Most often we will be using the contrapositive of the above definition.
\begin{remark} \label{R:projComp}
Let $C \subset \Hstar$ be a  projectively compact set with corresponding sets $C_n= \{ \vec{x}_n = (\ip<x,e_1>, \dots, \ip<x,e_n>) \, ; \; x\in C\}$. If  $x \notin C$,  then there exist an $N$ such that for $n > N$ we have that $x_n \notin C_n$.
\end{remark}

%
%

The following proposition discusses the properties of these projectively compact sets.

\begin{prop}
Let $C$ be a  projectively compact set with corresponding sets $C_n= \{ \vec{x}_n = (\ip<x,e_1>, \dots, \ip<x,e_n>) \, ; \; x\in C\} \subset \Rn$. Then $C$ is convex if and only if each $C_n$ is convex.
\end{prop}

\begin{proof}
Suppose $C$ is convex and consider the set $C_n$. Let $\vec{x}_n,\vec{y}_n$ be two points in $C_n$ corresponding to $x,y \in C$ . That is $\vec{x}_n = (\ip<x,e_1>, \dots, \ip<x,e_n>)$ and $\vec{y}_n = (\ip<y,e_1>, \dots, \ip<y,e_n>)$.  We must show $\ga \vec{x}_n + (1-\ga) \vec{y}_n$ is in $C_n$ for any $\ga \in [0,1]$. Since $C$ is convex we have that $\ga x + (1-\ga)y \in C$. Thus 
\[
(\ip<\ga x + (1-\ga)y ,e_1>, \dots, \ip<\ga x + (1-\ga)y , e_n>) \in C_n.
\]
Notice 
\[
(\ip<\ga x + (1-\ga)y ,e_1>, \dots, \ip<\ga x + (1-\ga)y , e_n>) =\ga \vec{x}_n + (1-\ga)\vec{y}_n
\]
and thus $\ga \vec{x}_n + (1-\ga)\vec{y}_n \in C_n$ and we have $C_n$ is convex. 

On the other hand suppose that $C_n$ is convex for each $n$. We must show that $C$ is convex. Let $x,y \in C$. We will show that $\ga x + (1-\ga)y \in C$. Since $x,y \in C$ we have that $\vec{x}_n, \vec{y}_n \in C_n$ for all $n$ and by the convexity of each $C_n$ we have that $\ga \vec{x}_n + (1-\ga) \vec{y}_n \in C_n$ for all $n$. Thus we have that $\ga x + (1-\ga)y \in C$.


\end{proof}

We will now demonstrate that there are nontrivial sets which satisfy this criteria of being  convex and {projectively  compact}. In particular we  demonstrate the  closed ball in $H_p'$ given by
\[
B^{-p}_r(y) = \{ x\in H_p' \, ; \, \norm|x-y|_{-p} \leq r \}
\]
has positive measure for any $r >0$ and $y \in H_p'$. Of course these sets are closed convex and  projectively  compact because their ``projections'' are essentially closed ellipses in $\Rn$. We must just demonstrate they have positive measure.

First we observe that every ball $B^{-p}_r(y)$ contains a ball centered around a ``rational point'' $q$, i.e. $q \in \text{span}_{\Q}\{e_1, \dots, e_n\} \subset \Hnucl$ for some $n$.

\begin{lemma} \label{L:balls}
For any $y \in \Hstar$ and $r >0$, the ball $B^{-p}_r(y)$ contains a ball $B^{-p}_{r'}(q)$ where $0 < r' < r$ and $q \in \text{span}_{\Q}\{e_1, \dots, e_n\}$ for some $n$. 
\end{lemma}

\begin{proof}
Consider the set $Q_n = \{r_1e_1 + \cdots r_ne_n \, ; \, r_1, \dots, r_n \in \Q \}$. Note that $Q_n$ is a countable set. Now let $Q = \bigcup_{n=1}^{\infty} Q_n$. Again $Q$ is countable. (This can be thought of as the set of rational points in $\Hstar$.) Observe if $y \in H_p'$, then there exists $n$ such that 
\[
\norm|y- \sum_{k=1}^{n} \ip<y,e_n>e_n|_{-p} < \frac{r}{2}
\]
For each $k =1, \dots, n$, take $r_k \in \Q$ such that 
\[
|r_k - \ip<y,e_k>| < \frac{r}{2^{k}}.
\]
Then let $q \in Q$ be given by $q = r_1e_1 + \cdots + r_ne_n$ and we have
\begin{align*}
\norm|y-q|_{-p} &\leq \norm|y-\sum_{k=1}^{n} \ip<y,e_k>e_k|_{-p} + \norm|\sum_{k=1}^{n} \ip<y,e_k>e_k-\sum_{k=1}^{n} r_ke_k|_{-p} \\
& < \frac{r}{2} + \frac{r}{2} = r
\end{align*}
Thus for any $x \in B^{-p}_{\frac{r}{2}}(q)$ we have 
\[
\norm|x-y|_{-p} \leq \norm|x-q|_{-p} + \norm|y-q|_{-p} \leq \frac{r}{2} + \frac{r}{2} \leq r.
\]
Therefore $B^{-p}_{\frac{r}{2}}(q) \subset B^{-p}_{r}(y)$
\end{proof}

We now use the previous lemma along with the properties of the measure $\mu$ to deduce that any such ball in $H_p'$ with positive radius must have positive measure.

\begin{prop} \label{P:positive}
For any $y \in H_p'$ and $r >0$ we have that $\mu(B_r^{-p}(y)) > 0$. 
\end{prop}

\begin{proof}
We let $Q$ be as in \lemref{L:balls}.
Since $Q$ is dense in $H_p'$ for any $r>0$, $H_p'$ can be written as a countable union of balls centered about rational points, i.e.
\[
H_p' = \bigcup_{q \in Q} B^{-p}_r(q).
\]
Since $H_p'$ is of positive measure (actually full measure), we have that $B^{-p}_r(q)$ must have positive measure for some $q$. By \propref{P:translation} we have that $B^{-p}_r(q)$ is of positive measure for any $q$. 

Every ball $B_r^{-p}(y)$ in $H_p'$ contains a ball centered at a rational point by \lemref{L:balls}. Thus each ball $B_r^{-p}(y)$ must have positive measure.
\end{proof}

The basis for  the inductive topology for $\Hstar$ is given by the convex hull of the sets
\[
\bigcup_{p=1}^{\infty} B^{-p}_{r_p}(x_p)
\] 
where $B^{-p}_{r_p}(x_p)$ denotes the unit ball in $H_p'$ centered around $x_p \in H_p'$ with radius $r_p$ \cite{JB06}.
Thus each nonempty open set contains an open ball $B^{-p}_r(x)$ for some $x \in H_{-p}'$ and $r >0$. Since each $B^{-p}_r(x)$ has positive measure by \propref{P:positive}, we must have that each nonempty open set in $\Hstar$ also has positive measure. 

Now the inductive and strong topologies on $\Hstar$ are equivalent and of course the weak topology is coarser than either of these topologies \cite{JB06}. This leads to the following corollary of \propref{P:positive}.

\begin{cor}
 The $\mu$-measure of any  nonempty open set $U$ in any of the weak, strong, or inductive limit topologies  is positive (i.e. $\mu(U) > 0$).
\end{cor}

 The next theorem is the main result of this paper. It gives us a Support Theorem for the Gauss--Radon Transform.

\begin{thm}[Support Theorem for Gauss--Radon Transform]
Let $\phi$  be a test function and let $C$ be a  convex projectively compact set in $\Hstar$ with
\[
G_{\phi}(\ga v + v^{\perp}) = 0
\]
when $\ga v \notin C$. Then $\phi(x) = 0$ for all $x \notin C$
\end{thm}

\begin{proof}
Throughout the proof we make the following notational conventions: for a vector  $w \in H_0$, we denote by $w_n$ the projection of $w$ onto $\text{span}\{e_1, \dots, e_n\}$. That is,
\[
w_n = \ip<w,e_1>e_1 + \ip<w,e_2>e_2 + \cdots + \ip<w, e_n>e_n.
\]
Also we denote by $\vec{w}_n$ the vector in $\Rn$ corresponding to $w_n$. That is,
\[
\vec{w}_n = (\ip<w,e_1>, \ip<w,e_2>, \cdots, \ip<w,e_n>).
\]
We also make the observation that if $w \notin C$, then there exist an $N$ such that  $w_n \notin C$ for all $n >N$. 
(For if $w_n \in C$ for all $n > N$, then $\vec{w}_n \in C_n$ for all $n$ and hence $w \in C$.)

With this in mind we 
first take a $\ga v \in V_n = \text{span}\{e_1, \dots, e_n\}$ with 
$\ga v \notin C$. We know that
\begin{equation}\label{E:GaussZero}
0=G_{\phi}(\ga v + v^{\perp})= \int_{\Hstar} \phi(x)  \,d\mu_{\ga v + v^{\perp}}(x).
\end{equation} 
We can decompose the subspace $v^{\perp}$ as follows
\[
v^{\perp} = (v^{\perp} \cap V_n) \oplus V_n^{\perp}.
\]
Let $S_{n,v} =(v^{\perp} \cap V_n) $ and notice $V_n^{\perp}=\text{span}\{e_{n+1}, e_{n+2}, \dots \}$. The idea here is that
$V_n$ is in essence $\Rn$ and thus basically we have $v$ ``in'' $\Rn$ and $S_{n,v}$ can be thought of as the orthogonal complement of $v$ ``in'' $\Rn$. 

We now use \thmref{T:dis} to rewrite $G_{\phi}(\ga v + v^{\perp})$ as follows:
\[
G_{\phi}(\ga v + v^{\perp}) = \int_{\Hstar}  \int_{\Hstar} \phi(x+y) \, d\mu_{V_n^{\perp}}(y) \, d\mu_{\ga v + S_{n,v}}(x).
\]
Let's write the inside of the above integral as 
\[
f_n(x) = \int_{\Hstar} \phi(x+y) \, d\mu_{V_n^{\perp}}(y)
\]
and we have that 
\begin{equation}\label{E:fn}
G_{\phi}(\ga v + v^{\perp}) = \int_{\Hstar} f_n(x) \, d\mu_{\ga v + S_{n,v}}(x).
\end{equation}
Since $v, S_{n,v}$ are in $\text{span}\{e_1, \dots, e_n\}$ the we can apply \thmref{T:coor} to write the above as
\begin{equation}\label{E:fnAgain}
\int_{\Hstar} f_n(x) \, d\mu_{\ga v + S_{n,v}}(x) = \int_{\Hstar} f_n(\ip<x,e_1>e_1 + \cdots + \ip<x,e_n>e_n) \, d\mu_{\ga v + S_{n,v}}(x).
\end{equation}
Moreover, $f_n(\ip<x,e_1>e_1 + \cdots + \ip<x,e_n>e_n) = F_n(\hat{e}_1, \dots, \hat{e}_n)$ where
\[
F_n(x_1, \dots, x_n) = \int_{\Hstar} \phi(x_1e_1 + \cdots + x_ne_n+y) \, d\mu_{V_n^{\perp}}(y)
\]
is a function on $\Rn$. Thus by \propref{P:Rn}
\[
\int_{\Hstar} f_n(\ip<x,e_1>e_1 + \cdots + \ip<x,e_n>e_n) \, d\mu_{\ga v + S_{n,v}}(x) = \int_{\ga \vec{v} + \vec{v}^{\perp}}  F_n \, d\mu_{\ga \vec{v} + \vec{v}^{\perp}}
\]
where $\vec{v} = (\ip<v,e_1>, \dots, \ip<v, e_n>)$ is a vector in $\Rn$ corresponding to $v$. Thus combining the above with \eqref{E:fn}  and \eqref{E:fnAgain} we have that
\[
 G_{\phi}(\ga v + v^{\perp})  =\int_{\ga \vec{v} + \vec{v}^{\perp}}  F_n \, d\mu_{\ga \vec{v} + \vec{v}^{\perp}}.
\] 
For any $\ga \vec{v} \notin C_n$ we must also that $\ga v \notin C$. Thus by assumption we have $ G_{\phi}(\ga v + v^{\perp}) =0$ and combining that with the above yields
\[
\int_{\ga \vec{v} + \vec{v}^{\perp}}  F_n \, d\mu_{\ga \vec{v} + \vec{v}^{\perp}} = 0.
\]
Therefore the Gauss--Radon transform of the function $F_n$ is $0$ when for any hyperplane $\ga \vec{v} + \vec{v}^{\perp}$ not intersecting the compact convex set $C_n$.

 We would like to apply the original Support Theorem (for the Gaussian measure) to our function $F_n$ to see that $F_n =0$ outside of $C_n$. In order to do so we must check that $F_n$ satisfies  the other assumptions of \thmref{T:STGaus}. 

\begin{lemma}
The function 
\[
F_n(x_1, \dots, x_n) = \int_{\Hstar} \phi(x_1e_1 + \cdots + x_ne_n+y) \, d\mu_{V_n^{\perp}}(y)
\]
is continuous and bounded.
\end{lemma}
\begin{proof}
The continuity is easy to check, observe that if $\{\vec{x}^{(k)} \}$ converges to $\vec{x}$ in $\Rn$, then
$x^{(k)} = x^{(k)}_1e_1 + \cdots + x^{(k)}_ne_n$ converges to $x = x_1e_1 + \cdots + x_ne_n$ with respect to $\norm|\cdot|_0$ (in fact, with respect to any $\norm|\cdot|_p$ or $\norm|\cdot|_{-p}$ norm). Therefore by \propref{P:deltaConv} we have that $\dV[x^{(k)} + V_n^{\perp}]$ converges to $\dV[x + V_n^{\perp}]$ strongly as $k \to \infty$. Therefore
\begin{align*}
\lim_{k \to \infty} F_n(x^{(k)}_1, \dots, x^{(k)}_n) &= \lim_{k \to \infty} \int_{\Hstar} \phi(x^{(k)}_1e_1 + \cdots + x^{(k)}_ne_n+y) \, d\mu_{V_n^{\perp}}(y) \\
&= \lim_{k \to \infty} \int_{\Hstar} \phi(y) \, d\mu_{x^{(k)}+V_n^{\perp}}(y) \qquad \text{by \thmref{T:shift}} \\
&= \lim_{k \to \infty} \pair<\phi, \dV[x^{(k)}+V_n^{\perp}]> \\
&= \pair<\phi, \dV[x+V_n^{\perp}]>  \qquad \text{by \propref{P:deltaConv} }\\
&= \int_{\Hstar} \phi(y) \, d\mu_{x+V_n^{\perp}}(y) \\
&=  \int_{\Hstar} \phi(x_1e_1 + \cdots + x_ne_n+y) \, d\mu_{V_n^{\perp}}(y) \qquad \text{by \thmref{T:shift}}
 \\
&= F_n(x_1,\dots, x_n)
\end{align*}
and thus $F_n$ is continuous on $\Rn$.

We now need to verify that $F_n$ is bounded. Again observe that
\[
F_n(x_1, \dots, x_n) = \int_{\Hstar} \phi(x_1e_1 + \cdots + x_ne_n+y) \, d\mu_{V_n^{\perp}}(y) = \int_{\Hstar} \phi(y) \, d\mu_{x+V_n^{\perp}}(y).
\]
Now $\phi$ is a test function and $\mu_{x + V_n^{\perp}}$ is a Hida measure. So we combining  \thmref{T:testBound} and \thmref{T:Hidabnd} to get that  for some $p \geq 1$ we have
\begin{align*}
|F_n(x_1, \dots, x_n)| &= \left| \int_{\Hstar} \phi(y) \, d\mu_{x+V_n^{\perp}}(y) \right| \\
  &\leq \int_{\Hstar} |\phi(y)| \, d\mu_{x+V_n^{\perp}}(y) \\ 
  &\leq  \int_{H_p'}  K_p \exp\left[\frac{1}{2} \norm|y|_{-p}^2 \right] d\mu_{x+V_n^{\perp}}(y)  
\end{align*}
and the last integral is finite by \thmref{T:Hidabnd}.
\end{proof}

Thus by the original Support Theorem (\thmref{T:STGaus})  we have that $F_n(\vec{x}) = 0$ for all 
$\vec{x} \notin C_n =\{(\ip<x,e_1>, \dots, \ip<x,e_n>) \, | \, x \in C \} \subset \Rn$. 
Thus
\[
f_n(x_n) = \int_{\Hstar} \phi(x_n+y) \, d\mu_{V_n^{\perp}}(y)
\]
is $0$ when  
$\vec{x}_n \notin C_n$.  To complete the proof we notice that $f_n(x_n) \to \phi(x)$. To this end notice that
\begin{align}
f_n(x_n) &= \int_{\Hstar} \phi(x_n+y) \, d\mu_{V_n^{\perp}}(y) \notag \\
         &=\int_{\Hstar} \phi(y) \, d\mu_{x_n +V_n^{\perp}}(y) \text{ by \thmref{T:shift}} \notag \\
         &=\pair<\phi, \dV[x_n + {V_n^{\perp}}]> \label{eq:fn}.
\end{align}

%
 Combining \propref{P:deltaConv} with \eqref{eq:fn} gives us that
\[
\lim_{n\to\infty} f_n(x_n) = \lim_{n\to\infty} \pair<\phi, \dV[x_n + {V_n^{\perp}}]> =\pair<\phi, \dV[x ]> = \phi(x).
\]
We now have the tools to complete the proof. Take an 
$x \notin C$. Then  by Remark \ref{R:projComp} 
 there exist an integer $N$ such that for all $n > N$,  we have  $\vec{x}_n \notin C_n$. Thus $F_n(\vec{x}_n) =0$
 and likewise $f_n(x_n) = 0$ for all $n > N$.
Taking the limit as $n$ goes to infinity gives us $\phi(x) = 0$.
\end{proof}

The following corollary restates the above when we simply take the closed ball $B_{r}^{-p}(x)$ as our closed convex projectively compact set.

\begin{cor}[Support Theorem for Gauss--Radon Transform]
Let $\phi$  be a test function with
\[
G_{\phi}(\ga v + v^{\perp}) = 0
\]
when $\norm|\ga v|_{-p} > r$. Then $\phi(x) = 0$ for all $x \notin B_{r}^{-p}(0)$
\end{cor}


\bibliographystyle{siam}
\bibliography{radon}  

\end{document}